\documentclass[12pt,twoside,draft,reqno]{amsart}

\usepackage{amssymb}

\newcommand{\ms}{\medskip}

\def\supp{\text{\rm supp}\,}
\def\to{\rightarrow}

\def\bs{\bigskip}\def
\ms{\medskip}

\def\OO{{\mathcal{O}}}

\def\R{{\mathbb R}}
\def\C{{\mathbb{C}}}

\def\HH{{\mathcal H}}

\def\e{\varepsilon}

\newtheorem{theorem}{Theorem}[section]
\newtheorem{lemma}[theorem]{Lemma}
\newtheorem{corollary}[theorem]{corollary}

\newtheorem{remark}[theorem]{Remark}

\title{Uniqueness theorems for Cauchy integrals}

\author{Mark~Melnikov}
\address{Department de Matematiques\\
Uneversitat Autonoma de Barcelona\\
08193 Bellaterra, Barcelona, Spain}
\email{mark.melnikov@gmail.com}
\thanks{The first
author is supported by grants No. MTM2004-00519
and 2001SGR00431}

\author{ Alexei~Poltoratski}
\address{Texas A\& M University
\\ Department of Mathematics\\
College Station, TX 77843, USA}
\email{alexeip@math.tamu.edu}
\thanks{The second
author is supported by 
N.S.F. Grant No. 0500852}

\author{Alexander~Volberg}
\address{Dept. Math. \\
Michigan State Univ.\\
East Lansing MI 48823, USA\\
and\\
School of Math.\\
University of Edinburgh\\
Edinburgh UK EH9 EJ6}
\email{sashavolberg@yahoo.com}
\thanks{The third
author is supported by 
N.S.F. Grant No. 0501067 }

\begin{document}
\maketitle
\begin{abstract}
If $\mu$ is a finite complex measure in the complex plane $\C$ we denote by $C^\mu$ its Cauchy integral
defined in the sense of principal value. The measure $\mu$ is called {\it reflectionless} if it is continuous (has no atoms) and $C^\mu=0$ at $\mu$-almost every point.
We show that if $\mu$ is reflectionless and its Cauchy maximal function $C^\mu_*$ is summable with respect to $|\mu|$
then $\mu$ is trivial. An example of  a reflectionless measure whose maximal function belongs to the "weak"
$L^1$ is also constructed, proving that the above result is sharp in its scale. We also give a partial geometric description
of the set of reflectionless measures on the line and discuss connections of our results with the notion of sets of finite
perimeter in the sense of De Giorgi.
\end{abstract}
\section{Introduction}

This article discusses uniqueness theorems for Cauchy integrals of complex measures in the plane. 
We consider the space $M=M(\C)$ of finite complex measures $\mu$ in $\C$.
The Cauchy integral of a measure from $M$ is defined in the sense of principal value. First, for any  $\mu\in M$, $\varepsilon>0$ and any $z\in \C$  consider 
$$C^{\mu}_{\varepsilon}(z):=\int_{\zeta: |\zeta-z|>\varepsilon}\frac{d\mu(\zeta)}{\zeta-z}.$$
Consequently, the Cauchy integral of $\mu$ can be defined as 
$$
C^{\mu}(z):= \lim_{\varepsilon\rightarrow 0}C^{\mu}_{\varepsilon}(z) \,,
$$
if the limit exists. 

Unlike the Cauchy transform on the line, $C^\mu$ can vanish on a set of positive Lebesgue measure: consider for example $\mu=dz$ on a closed curve,
whose Cauchy transform is zero at all points outside the curve. It is natural to ask if $C^\mu$ can also vanish on large sets with respect to $\mu$. 
 If $\mu=\delta_z$ is a single point mass, its Cauchy transform will be zero $\mu$-a.e. due to the above definition
of $C^\mu$ in the sense of principal value. Examples of infinite discrete measures with vanishing Cauchy transforms can also be constructed with little
effort. 

After that one arrives at the following corrected version of the question: Is it true that any continuous $\mu\in M$, such that $C^{\mu}(z)=0$ at $\mu$-a.e.  point, is trivial? 
As usual, we call a measure continuous if it has no point masses. We denote the space of all finite complex continuous measures by $M_c(\C)$.

This problem can also be interpreted in terms of uniqueness. Namely, if $f$ and $g$ are two
functions from $L^1(|\mu|)$ such that $C^{(f-g)\mu}=0$, $\mu$-a.e., does it imply that $f= g$, $\mu$-a.e.? This way it becomes a problem of injectivity of the planar
Cauchy transform.

First significant progress towards the solution of this problem was achieved by X. Tolsa and J. Verdera in \cite{TV}. It was established that the answer is positive in two important particular cases:
when $\mu$ is absolutely continuous with respect to Lebesgue measure $m_2$ in $\C$ and when $\mu$ is a measure of linear growth with finite Menger curvature. The latter class of measures is one of  the main objects in the study of the planar Cauchy transform, see for instance \cite{MV}, \cite{NTV} or \cite{T}.

As to the complete solution to the problem, it seemed for a while that the answer could be positive for any  $\mu\in M_c$, see for example \cite{TV}.
However, in Section \ref{refl} of the present paper we show that there exists a large set of  continuous measures $\mu$ satisfying $C^{\mu}(z)=0$, $\mu$-a.e.
Following \cite{BBEIM}, we call such measures {\it reflectionless}. 
This class seems to be an intriguing new object in the theory.

On the positive side, we prove  that if the maximal function associated with the Cauchy transform is summable with respect to 
$|\mu|$
then $\mu$ cannot be reflectionless, see Theorem \ref{zero}. This result is sharp in its scale because the simplest examples of reflectionless measures produce maximal functions that lie in 
the "weak" $L^1(|\mu|)$.
We prove this result in Section \ref{1}
 In view of this fact, we believe that the class of continuous measures with summable Cauchy maximal functions also deserves
attention.

\ms

 A full description of this class and the (disjoint) class of reflectionless measures remains an {\bf open problem}.

\ms

Let us mention that if $\mu$ is a measure with linear growth and finite Menger curvature then its Cauchy maximal function belongs to $L^2(|\mu|)$, see \cite{NTV, T}, and
therefore is summable.
This fact relates Theorem \ref{zero} to the beforementioned result from \cite{TV}. The latter can also be deduced in a different way, see Section \ref{1}.

From the point of view of uniqueness, our results imply that any bounded planar Cauchy transform is injective, see corollary \ref{c2}. This property is a clear analogue of
the uniqueness results for the Cauchy integral on the line or  the unit circle.

In Section \ref{DG} we discuss other applications of Theorem \ref{meln}. They involve structural theorems of De Giorgi and his notion of a set of finite perimeter,
see \cite{EG}.

 In Section \ref{2} we study asymptotic behavior of the Cauchy transform near its zero set. The results of this section imply that the Radon derivative of $\mu$  with respect to Lebesgue measure $m_2$ vanishes a.e. on the set $\{C^\mu=0\}$. In particular the set  $\{C^\mu=0\}$ must be a zero set with respect to the variation of the
absolutely continuous part of $\mu$ which is a slight generalization of the first result of \cite{TV}.  It is interesting to note that the most direct analogue of this corollary on the real line is false: it is easy to construct an absolutely  continuous (with respect to $m_1=dx$) measure $\mu\in M(\R)$ such that $|\mu|(\{C^\mu=0\})>0$.

Finally, in Section \ref{refl} we attempt a geometric description of the set of reflectionless measures. We give a partial description
of reflectionless measures on the line in terms of so-called comb-like domains. We also provide  tools for the construction of various examples of such measures.
In particular, we show that the harmonic measure on any compact subset (of positive Lebesgue measure) of $\R$ is reflectionless.

\ms\noindent
{\bf Acknowledgments.} The authors are grateful to Fedja Nazarov for his invaluable comments and insights. The second author would also like to thank
the administration and staff of Centre de Recerca Matem\'atica in Barcelona for the hospitality during his visit in the Spring of 2006. 
\vspace{.2in}

\section{Measures with summable maximal functions}\label{1}

If $\mu\in M$ we denote by $C^{\mu}_*(z)$ its Cauchy maximal function
$$
C^{\mu}_*(z):= \sup_{\varepsilon> 0} |C^{\mu}_{\varepsilon}(z)|.
$$

Our first result is the following uniqueness theorem. 

\begin{theorem}
\label{zero} Let $\mu\in M_c$.
Assume that  $C^{\mu}_*(z)\in L^1(|\mu|)$ and that  $C^{\mu}(z)$ exists and vanishes  
$\mu$-a.e. Then $\mu\equiv 0$.
\end{theorem}
\vspace{.2in}

We first prove

\vspace{.1in}

\begin{theorem}
\label{meln}
If $C^{\mu}_* \in L^1(|\mu|)$ and  $C^{\mu}(z)$ exists $\mu$-a.e. 
then
\begin{equation}
\label{quadratic}
2C^{C^{\mu}d\mu}(z) = 2\int\frac{ C^{\mu}(t)d\mu(t)} {t-z}=\left[C^\mu(z)\right]^2 \,\,\,\text{for}\,\,\, m_2\text{-a.e. point}\,\,z\in \C\,.
\end{equation}
\end{theorem}

\vspace{.1in}

\begin{proof} Put
$$
F:=  \{z\in \C: \int\frac{d|\mu|(t)}{|t-z|} <\infty\}\,.
$$
As $|\mu|$ is a finite measure,
\begin{equation}
\label{3}
m_2(\C \setminus F) =0\,.
\end{equation}
Let $z\in F$.
Then  the integral
$$
I:=\int\int_{|t-\zeta|>\varepsilon} d\mu(t)d\mu(\zeta)\frac1{t-z}\cdot\frac1{\zeta-z}
$$
is absolutely convergent for any $\varepsilon>0$.

Using the identity
$$\frac1{(t-z)(z-\zeta)}+\frac1{(z-\zeta)(\zeta-t)}+\frac1{(\zeta-t)(t-z)}\equiv 0$$
we obtain
$$
I = \int\int_{ |t-\zeta|>\varepsilon}\bigg[\frac1{z-\zeta}\cdot\frac1{\zeta-t}+\frac1{\zeta-t}\cdot\frac1{t-z}\bigg]d\mu(t)d\mu(\zeta)=
$$
$$
\int\frac{d\mu(\zeta)}{\zeta-z}\int_{|t-\zeta|>\varepsilon}\frac{d\mu(t)}{t-\zeta}+
\int\frac{d\mu(t)}{t-z}\int_{ |\zeta-t|>\varepsilon}\frac{d\mu(\zeta)}{\zeta-t} 
 =
$$
$$
\int d\mu(t)\cdot C^{\mu}_{\varepsilon}(t)\cdot \frac1{t-z} + \int d\mu(\zeta)\cdot C^{\mu}_{\varepsilon}(\zeta)\cdot \frac1{\zeta-z} = 2\int \frac{ C^{\mu}_{\varepsilon}(t)d\mu(t)} {t-z}\,.
$$
Put
$$
E:=  \{z\in \C: \int\frac{C^{\mu}_{*}(t)d|\mu|(t)}{|t-z|} <\infty\}\,.
$$
By assumption, the numerator $C^{\mu}_*(t)d|\mu|(t)$ is a finite measure. Therefore
\begin{equation}
\label{e1}
m_2(\C \setminus E) =0\,.
\end{equation}

If $z\in E$ then
\begin{equation}
\label{e20}
\lim_{\varepsilon\rightarrow 0}\int\frac{ C^{\mu}_{\varepsilon}(t)d\mu(t)} {t-z} =\int \frac{ C^{\mu}(t)d\mu(t)} {t-z}\,.
\end{equation}

This formula is true as long as $C^{\mu}_* \in L^1(|\mu|)$ and the principal value $C^{\mu}$ exists $\mu$-a.e. 
by the dominated convergence theorem. Thus
\begin{equation}
\label{2a}
\lim_{\varepsilon\rightarrow 0}I = 2C^{C^{\mu}d\mu}(z)\,\,\,\text{if}\,\,\,\,z\in E \,.
\end{equation}

It is left to show that, since $z\in F$,
\begin{equation}
\label{4}
\lim_{\varepsilon\rightarrow 0}I  =[C^{\mu}(z)]^2\,.
\end{equation}
Since $z\in F$,  the following integral converges absolutely:
$$
\phi_{\varepsilon}(t,z) := \int_{\zeta\in \C, |\zeta-t|>\varepsilon} \frac{d\mu(\zeta)}{\zeta-z}\,.
$$
Also
$$
I = \int \phi_{\varepsilon} (t,z) \frac1{t-z} d\mu(t)\,.
$$
Since the point $z$ is fixed in $F$, we have that $\frac1{|\zeta-z|}\in L^1(|\mu|)$, and therefore
$\int_A\frac1{|\zeta-z|}d|\mu|(\zeta)$ is small if $|\mu|(A)$ is small. Denoting the disc centered at $t$ and of radius $\varepsilon$ by $B(t,\varepsilon)$ we notice that
$$
1)\,\, \phi_{\varepsilon}(t,z) = \int_\C \frac{d\mu(\zeta)}{\zeta-z} - \int_{B(t,\varepsilon)} \frac{d\mu(\zeta)}{\zeta-z}\,,
$$
\bs
$$
2)\,\, \lim_{\varepsilon\rightarrow 0} |\mu|(B(t,\varepsilon)) =0.
$$
uniformly in $t$. Otherwise $\mu$ would have an atom.

We conclude that, as $\e\to 0$, the functions
$\phi_{\varepsilon}(t,z)$ converge uniformly in $t\in \C$ to $\phi(z)=\int \frac{d\mu(\zeta}{\zeta-z}$. Hence for any $z\in F$ and any $t\in \C\setminus{z}$
$$
3)\,\, \frac{\phi_{\varepsilon}(t,z)}{t-z} \rightarrow \frac{\phi(z)}{t-z}, \,\, \text{as}\ \ \varepsilon\rightarrow 0\,.
$$
Since $\phi_{\varepsilon}(t,z)$ converge uniformly and $z\in F$, 
$$
\int d\mu(t) \phi_{\varepsilon}(t,z) \frac1{t-z} \rightarrow \phi(z)\int\frac{d\mu(t)}{t-z}= [C^{\mu}(z)]^2\,.
$$
We have verified \eqref{4}.

\vspace{.1in}

Combining \eqref{2a} and \eqref{4} we conclude that
for $z\in E\cap F$ (so for $m_2$-a.e. $z\in \C$) we have
\begin{equation}
\label{quadratic0}
2C^{C^{\mu}d\mu}(z) = 2\int\frac{ C^{\mu}(t)d\mu(t)} {t-z} =\lim_{\varepsilon\rightarrow 0} I =[C^{\mu}(z)]^2\,\,\,\text{for}\,\,\, m_2\text{-a.e. point}\,\,z\in \C\,.
\end{equation}

\vspace{.1in}

This formula is true as long as $C^{\mu}_* \in L^1(|\mu|)$ and the principal value $C^{\mu}$ exists $\mu$-a.e. 

\end{proof}

\vspace{.2in}

To deduce Theorem \ref{zero}
suppose that $C^\mu$ vanishes $\mu$-a.e. Then the left-hand side in \eqref{quadratic0} is zero for $m_2$-a.e. point $z$. The same must hold for $[C^{\mu}(z)]^2$.
But if $C^{\mu}(z) =0$ for Lebesgue-a.e. point $z\in \C$ then $\mu=0$, see for example \cite{Gam}. Theorem \ref{zero} is completely proved.

\bs

{\bf Remark.} In the statement of Theorem \ref{meln} the condition $C^\mu_*\in L^1(|\mu|)$ can be replaced with the condition that $C_\e^\mu$
converge in $L^1(|\mu|)$. The proof would have to be changed as follows. 

Like in the above proof one can show that at Lebesgue-a.e. point $z$
\begin{equation}
\lim_{\varepsilon\rightarrow 0}I  =[C^{\mu}(z)]^2\,.
\end{equation}
The relation 
$$I= 2\int \frac{ C^{\mu}_{\varepsilon}(t)d\mu(t)} {t-z}
$$
for a.e. $z$ can also be established as before. Since $C_\e^\mu$
converge in $L^1(|\mu|)$, the last integral converges to $C^{C^\mu d\mu}(z)$ in the "weak" $L^2(dxdy)$, which concludes the proof.

\ms

Hence we arrive at the following version of Theorem \ref{zero}:

\ms

\begin{theorem}\label{zerom}
Let $\mu\in M_c$. Assume that  $C^{\mu}_\e\to 0$ in $ L^1(|\mu|)$. Then $\mu\equiv 0$.
\end{theorem}

This version has the following corollary:

\begin{corollary}[\cite{TV}]\label{c1}
Let $\mu\in M$ be a measure of linear growth and finite Menger curvature. If $C^\mu=0$ at $\mu$-a.e. point then $\mu\equiv 0$.
\end{corollary}

\begin{proof} The conditions on $\mu$ imply that the $L^2(|\mu|)$-norms of the functions $C^\mu_\e$ are uniformly bounded, see for instance \cite{MV}.
Since  $C^\mu_\e$ also converge $\mu$-a.e., they must converge in $L^1(|\mu|)$.
\end{proof}

{\bf Remark} As was mentioned in the introduction, Corollary \ref{c1} also follows from Theorem \ref{zero}. However, the above version of the argument allows one
to obtain it without the additional results of \cite{NTV, T} on the maximal function.

\bs

We also obtain the following statement on the injectivity of any bounded planar Cauchy transform. As usual, we say
that the Cauchy transform is bounded in $L^2(\mu)$ if the functions $C^{fd\mu}_\e$ are uniformly bounded
in $L^2(\mu)$-norm for any $f\in L^2(\mu)$. If $C^\mu$ is bounded, then
$C^{fd\mu}_\e$ converge $\mu$-a.e as $\e\to 0$ and the image $C^{fd\mu}$ exists in a regular sense
as a function in $L^2(\mu)$, see \cite{T}.

\begin{corollary}\label{c2} Let $\mu\in M$ be a positive measure. If $C^\mu$ is bounded in $L^2(\mu)$ then it is injective (has a trivial kernel).
\end{corollary}

\begin{proof} Suppose that there is $f\in L^2(\mu)$ such that  $C^{fd\mu}=0$ at $\mu$-a.e. point. Since both $f$ and $C^{fd\mu}_*$ are in $L^2(\mu)$,
$C^{fd\mu}_*$ is in $L^1(|f|d\mu)$.
Hence $f$ is a zero-function by Theorem \ref{zero}
\end{proof}

{\bf Remark} We have actually obtained a slightly stronger statement:  If $C^\mu$ is bounded in $L^2(\mu)$ then for any $f\in L^2(\mu)$
the functions $f$ and $C^{fd\mu}$ cannot have disjoint essential supports, i.e. the product $fC^{fd\mu}$ cannot equal to 0 at $\mu$-a.e. point.

\bs\bs

In the rest of this section we will discuss what other kernels could replace the Cauchy kernel in the statement of Theorem \ref{zero}.

If $K(x)$ is a complex-valued function in $\R^n$, bounded outside of any neighborhood of the origin, and $\mu$ is a finite measure on $\R^n$,
one can define $K^\mu$ and $K^\mu_*$ in the same way as $C^\mu$ and $C^\mu_*$ were defined in the introduction.

The proof of Theorem \ref{meln} relied on  the fact that the Cauchy kernel $K(z)=1/z$ is odd, satisfies the symmetry condition \eqref{e1}, i.e. 
\begin{equation}K(x-y)K(y-z)+K(y-z)K(z-x)+K(z-x)K(x-y)\equiv 0,\label{e2}\end{equation}
and is summable as
a function of $z$ for any $t$ with respect to Lebesgue measure. Any $K(x)$ having these three properties could be used in Theorem \ref{zero}.
Out of these three conditions the symmetry condition \eqref{e2} seems to be most unique. However, other symmetry conditions may result
in formulas similar to Theorem \ref{meln} that could still yield Theorem \ref{zero}. 

Here is a different example. It shows that much less symmetry can be required from the kernel if the measure is positive. 

\begin{theorem}\label{t3}
Let $\mu$ be a positive measure in $\R^n$. Suppose that the real kernel $K(x)$ satisfies the following properties:

1) $ K(-x)=- K(x)$ for any $x\in \R^n$;

2) $ K(x)>0$ for any $x$ from the half-space $\R^n_+=\{x=(x_1,x_2,...,x_n)\ |\ x_1>0\}$.

If $ K^\mu_*\in L^1(\mu) $ and $ K^\mu(x)=0$ for $\mu$-a.e. $x$ then $\mu\equiv 0$.

\end{theorem}

Note that  real and imaginary parts of the Cauchy kernel, Riesz kernels in $R^n$, as well as many other standard kernels satisfy the conditions
of the theorem.

We will need the following

\begin{lemma} Let $K$ be an odd kernel. 
and let $\mu,\nu\in M$. Then
\begin{equation}\int K^\mu_\e (z) d\nu(z)=-\int K^\nu_\e (z) d\mu(z)\end{equation}
for any $\varepsilon>0$.

Suppose that $K_*^\mu\in L^1(|\nu|)$. If $K^\mu(z)$ exists $\nu$-a.e.
then

$$\int K^\mu (z) d\nu(z)=-\lim_{\e\to 0}\int K^\nu_\e (z) d\mu(z).$$

In particular, suppose that both $K_*^\mu\in L^1(|\nu|)$ and $K_*^\nu\in L^1(|\mu|)$. If 
$K^\mu(z)$ exists $\nu$-a.e. and $K^\nu(z)$ exists $\mu$-a.e. then

$$\int K^\mu (z) d\nu(z)=\int K^\nu (z) d\mu(z).$$

\end{lemma}

\begin{proof} Since $K$ is odd, the first equation can be obtained simply by changing the order of integration.
 The second and third equations now follow from the dominated convergence theorem.
\end{proof}

\begin{proof}[Proof of Theorem \ref{t3}]
There exists a half-plane $\{x_1=c\}$ in $\R^n$ such that 
$\mu(\{x_1=c\})=0$ but both
$\mu(\{x_1>c\})$ and $\mu(\{x_1<c\})$ are non-zero.  Denote by $\nu$ and $\eta$ the restrictions of $\mu$ onto $\{x_1>c\}$
and $\{x_1<c\}$ respectively.
Then 
$$\int K^\nu_\e(z)d\mu(z)=\int K^\nu_\e(z)d\nu(z)+\int K^\nu_\e(z)d\eta(z).$$
The first integral on the right-hand side is 0 because of the oddness of $ K$ (apply the first equation in the last lemma with $\mu=\nu$). The second condition on $K$ 
and the positivity of the measure
imply that the second integral is positive and increases as $\e\to 0$. Therefore $\int K^\nu_\e(z)d\mu(z)$ cannot tend to zero. This contradicts  the fact that $K^\mu=0$,
$\nu$-a.e. and the second equation from the last lemma.
\end{proof}

\section{ Sets of finite perimeter }\label{DG}


In this section we give another example of an application of Theorem \ref{meln}. It involves the notion of  a set of finite perimeter introduced by De Giorgi in the 50's, see \cite{EG}. 
We say that a set $G\subset\R^2$ has finite perimeter (in the sense of De Giorgi) if the distributional partial derivatives of its characteristic function $\chi_G$ are finite measures.
Such sets have structural theorems. For example, if $G$ is such a set then the measure $\nabla\chi_G$ is carried by a set $E$, rectifiable in the sense of Besicovitch, i. e. a subset of a countable union of $C^1$ curves and an $\HH^1$-null set, where $\HH^1$ is the one-dimensional Hausdorff measure. Also the measure $\nabla\chi_G$ is absolutely continuous with respect to $\HH^1$ restricted to $E$ and its Radon-Nikodym derivative
is a unit normal vector $\HH^1$-a.e. (notice that $\nabla\chi_G$ is a vector measure). At $\HH^1$-almost all points of $E$ the function $\chi_G$ has approximate ``one-sided"' limit. For more details we refer the reader to  \cite{EG}.

The general question we consider can be formulated as follows: What can be said about $\mu$ if $C^\mu$ coincides at $\mu$-a.e. point with a "good" function $f$? To avoid certain technical details, all measures  in this section are compactly supported. 
Furthermore,  we will only discuss the two simplest choices of $f$. As we will see, even in such elementary situations
Theorem \ref{meln} yields interesting consequences.

As usual, when we say that $C^\mu=f$ at $\mu$-a.e. point, we imply that the principal value exists $\mu$-almost everywhere. 

\vspace{.1in}

\begin{theorem}
\label{one}
Let  $\mu\in M_c$ be compactly supported. Assume   that $C^{\mu}(z)=1$, 
$\mu$-almost everywhere  and $C^\mu_*\in L^1(|\mu|)$. Then $\mu=\bar{\partial}\chi_G$, where $G$ is a set of finite perimeter. In particular, $\mu$ is carried by a set $E$, $\HH^1(E) <\infty$, rectifiable in the sense of Besicovitch, and $\mu$ is absolutely continuous with respect to the restriction of $\HH^1$ to $E$. 
\end{theorem}

\vspace{.1in}

\noindent{\bf Remark.} The most natural example of such a measure is $dz$ on a $C^1$ closed curve. The theorem says that, by the structural results of De Giorgi, this is basically the full answer.

\vspace{.1in}

\begin{proof}

By Theorem \ref{meln} we get that for Lebesgue-almost every point in $\C$ 
\begin{equation}
\label{quadratic1}
[C^{\mu}(z)]^2 = 2\,C^{\mu}(z)\,.
\end{equation}

\vspace{.1in}

In other words for $m_2$-a.e. point $z$ we have $C^{\mu}(z)=0$ or $=2$. Let $G$ denote the set where $C^{\mu}(z) =2$. Since the Cauchy transform of any
compactly supported finite measure must tend to zero at infinity, this set is bounded. Consider the following equality 
$$
\chi_G=C^{\mu/2} ,
$$
understood in the sense that the two functions are equal as distributions.
Taking distributional derivatives on both sides we obtain
$$
\bar{\partial}\chi_G = \mu/2\ \  \text{and}\ \  \partial\bar\chi_G = \bar\mu/2.
$$
Hence $G$ has finite perimeter and the rest of the statement follows from the results of \cite{EG}.
\end{proof}

\bs

We say that a set $G$ has locally finite perimeter (in the sense of De Giorgi) if the distributional derivatives
of $\chi_G$ are locally finite measures.
Our second application is the following

\begin{theorem}
\label{notzero}
Let $\mu\in M_c$ be compactly supported. Assume   that $C^{\mu}(z)=z$, 
$\mu$-almost everywhere  and $C^\mu_*\in L^1(|\mu|)$. If $\mu(\C)=0$ then
$\mu=2z\bar{\partial}\chi_G$, where $G$ is a set with locally finite perimeter. Whether $\mu(\C)=0$ or not, $\mu$ is carried by a set $E$, $\HH^1(E) <\infty$, which is a rectifiable set in the sense of Besicovitch, and $\mu$ is  absolutely continuous with respect to the restriction of $\HH^1$ to $E$. 
\end{theorem}

\vspace{.1in}

\noindent{\bf Remark.} The most natural example of such a measure is $zdz$ on a $C^1$ closed curve. Our statement shows that this is basically one-half
of the  answer. The other half is given by $\sqrt{z^2-c} dz$ as will be seen from the proof.

\vspace{.1in}

\begin{proof}
Again, from Theorem \ref{meln} we get that for Lebesgue-almost every point in $\C$ 
\begin{equation}
\label{quadratic2}
[C^{\mu}(z)]^2 = 2\,C^{\zeta d\mu(\zeta)}(z)\,.
\end{equation}

Notice that
$$
C^{\zeta d\mu(\zeta)}(z) =\int \frac{\zeta}{\zeta-z}d\mu(\zeta)=\mu(\C) + z C^{\mu}(z)
$$
and we get a quadratic equation
$$
[C^{\mu}(z)]^2 = 2z C^{\mu}(z) -p\,,
$$ 
where $p:= -2\mu(\C)$.

\vspace{.1in}

First case $p=0$.
Here  we get
$$
[C^{\mu}(z)]^2 = 2z C^{\mu}(z)\,.
$$ 

We conclude that $C^{\mu}(z)= 0$ or $z$ for Lebesgue-a.e. point $z\in \C$.

Again a bounded set $G$ appears on which
$$
C^{\mu} = 2z\chi_G(z)
$$
in terms of distributions.
Therefore
$$
\bar{\partial}\chi_G =d\mu/2z\,,
$$
and the right hand side is a finite measure on any compact set avoiding the origin. Therefore, $G$ is a (locally)  De Giorgi set.

\vspace{.2in}

Let us consider the case $p\neq 0$. For simplicity we assume $p=1$, other $p$'s are treated in the same way. Then we have to solve the quadratic equation 
$$
C^\mu(z)^2 -2z C^\mu(z) +1=0
$$
for Lebesgue-a.e. point in $\C$. 
Let us make the slit $[-1,1]$ and consider two holomorphic functions in $\C\setminus [-1,1]$
$$
r_1(z)=z-\sqrt{z^2-1},\,\, r_2(z)=z+\sqrt{z^2-1}\,,
$$
where the branch of the square root is chosen so that 
$$
r_1(z)\rightarrow 0,\,\, z\rightarrow\infty\,.
$$
In other words we have the sets $E_1$ and $E_2$ such that
$m_2(\C\setminus E_1\cup E_2) =0$ and 
$$
z\in E_1\Rightarrow C^{\mu}(z) = r_1(z)\,,
$$
$$
z\in E_2\Rightarrow C^{\mu}(z) =  r_2(z)\,.
$$

\vspace{.1in}

Obviously it is $E_1$ that contains a neighborhood of infinity.
The function $z-\sqrt{z^2-1}$ outside of $[-1,1]$ can be written as $C^{\mu_0}(z)$ where $d\mu_0(x) = \frac1{\pi }\sqrt{1-x^2}dx$.
Consider $\nu=\mu-\mu_0$. Then 
$$
z\in E_1\Rightarrow C^{\nu}(z) = 0\,,
$$

\vspace{.1in}

$$
z\in E_2\Rightarrow C^{\nu}(z) = 2\sqrt{z^2-1}:=R(z)\,.
$$
Therefore,
\begin{equation}
\label{R}
C^{\nu}(z) =R(z)\chi_{E_2}\,.
\end{equation}
Notice that if $R$ was analytic  in an open domain compactly containing $E_2$ we would conclude from the previous equality that
$$
\nu = R(z) \bar{\partial}\chi_{E_2}.
$$
If, in addition, $|R|$ was bounded away from zero on $E_2$, we would obtain that
$\bar{\partial}\chi_{E_2}$ and $\partial\chi_{E_2}$ are measures of finite variation, and hence $E_2$ is a set of finite perimeter. Notice that our $R(z) =2\sqrt{z^2-1}$ is analytic in $O:=\C\setminus [-1,1]$ and is nowhere zero. We will conclude that $E_2$ is a set of locally finite perimeter.
More precisely we will establish the following claim:

\vspace{.1in}

\noindent For every open disk $V\subset O$ the set $O\cap E_2$ has finite perimeter.

\vspace{.1in}

Indeed, let $W$ be a disk compactly containing $V$, $W\subset O$. Let $\psi$ be a smooth function, supported in $W$, $\psi|V=1$. Multiply \eqref{R} by $\psi$ and take a distributional derivative (against smooth functions supported in $V$). Then we get (using the fact that $R$ is holomorphic on $V$)
$$
\nu|V = \bar{\partial}(\psi R \chi_{E_2\cap V})|V=\bar{\partial}( R \chi_{E_2\cap V})|V=R\bar{\partial}(\chi_{E_2\cap V})|V\,.
$$
We conclude immediately that $E_2\cap V$ is a set of finite perimeter. Therefore, $E_2\cap D$ is a set of finite perimeter, where $D$ is a domain whose closure is contained compactly in $O$.

\vspace{.1in}

Recalling that $\mu=\nu +\mu_0$ we finish the proof. \end{proof}

\vspace{.2in}

\begin{remark}
In is interesting to note that, as follows from the proof, if $\mu$ is the measure from the statement of the theorem then one of the connected components
of $\supp\mu$ must contain both roots of the equation $z^2+2\mu(\C)=0$.
\end{remark}

\vspace{.2in}

We conclude this section with the following examples of measures $\mu$ whose Cauchy transform coincides with $z$ at $\mu$-a.e. point

\bs

\noindent{\bf Examples}. 
\noindent 1. Let $\Omega$ be an open domain with smooth boundary $\Gamma$. Suppose that $[-1,1]\subset\Omega$. Let $\{D_j\}_{j=1}^{\infty}$ be smoothly bounded disjoint domains in $\OO:=\Omega\setminus [-1,1]$, $\gamma_j =\partial D_j$. Assume
\begin{equation}
\label{finitesum}
\sum_j \HH^1(\gamma_j) <\infty\,.
\end{equation}
Let $R(z)$ be  an analytic branch of $2\sqrt{z^2-1}$ in $\OO$. Consider the measure $\nu$  on
$\Gamma\cup(\cup\gamma_j)\cup [-1,1]$ defined as 
$$\nu= R(z) dz|_\Gamma-R(z) dz|_{\cup\gamma_j}-\frac1{\pi}\sqrt{1-x^2}dx|_{[-1,1]}.$$ Then 
$$
C^{\nu}(z)=
\begin{cases} &0\,\,\text{if}\,\, z\in \C\setminus \bar{\OO}\,,\\
&0\,\,\text{if}\,\,z\in\cup_j D_j\,,\\
&R(z)\,\,\text{if}\,\,z\in \OO\setminus\cup_j \bar{D}_j\,.
\end{cases}
$$
Recall that $R(z) = z+\sqrt{z^2-1}-(z-\sqrt{z^2-1})$ and that $C^{\mu_0}(z)=z-\sqrt{z^2-1}$ for $\mu_0=\frac1{\pi}\sqrt{1-x^2}dx|_{[-1,1]}$. We conclude that
for $\mu=\nu+\mu_0$ one has
$$
C^{\mu}(z)=
\begin{cases} & z-\sqrt{z^2-1}\,\,\text{if}\,\, z\in \C\setminus \bar{\OO}\,,\\
&z-\sqrt{z^2-1}\,\,\text{if}\,\,z\in\cup_j D_j\,,\\
&z+\sqrt{z^2-1}\,\,\text{if}\,\,z\in \OO\setminus\cup_j \bar{D}_j\,.
\end{cases}
$$

\vspace{.1in}

\noindent 2. The second example is exactly the same as the first one but $D_{j,k}= B(x_{j,k}, \frac1{10j^2})$, $x_{j,k}=2+\frac1{j}e^{\frac{2\pi i k}{j}}$, $1\le k\le j$,
$j=1,2,3...$. Here the assumption \eqref{finitesum} fails. But  $\nu$, defined as above, will still be a measure of finite variation (and so will be $\mu$): $|\nu|(\C) \leq C\sum_j\frac1{j^{3/2}}$.

\vspace{.1in}

In both examples $C^{\mu}(z) =z$ for $\mu$-a.e. $z$.

\section{Asymptotic behavior near the zero-set of $C^\mu$}
\label{2}

In this section we take a slightly different approach. We study asymptotic properties
of measures near the sets where the Cauchy transform vanishes. Theorem \ref{t4} below shows that near the density points of such sets the measure
must display a certain "irregular" asymptotic behavior.

As was mentioned in the introduction, one of the results of \cite{TV} says that
an absolutely continuous planar measure cannot
be reflectionless. This result is not implied by our Theorem \ref{zero} because
an absolutely continuous measure may not have a summable
Cauchy maximal function. It is, however, implied by Theorem \ref{t4}, see Corollary \ref{V}
below.

When estimating Cauchy integrals one often uses an elementary observation that  the difference of any two  Cauchy kernels
$1/(z-a) - 1/(z-b)$ can be estimated as $O(|z|^{-2})$ near infinity. To obtain higher order of 
decay one may consider higher order differences. Here we will utilize the following estimate of that kind, which can be
verified through simple calculations.

\begin{lemma}
If $a,b,c\in B(0,r)$ be different points, $|a-b|>r$.  Then there exist constants $A,B\in \C$ such that $|A|,|B|<2$

\begin{equation}\left|\frac A{z-a}+\frac B{z-b}-\frac 1{z-c}\right|<\frac {Cr^2}{|z|^3}\label{est}\end{equation}

outside of $B(0,2r)$. 

(Namely, $A=\frac{b-c}{b-a}, B=\frac{a-c}{a-b}$.)
\end{lemma}

If $\mu\in M$ consider one of its Riesz transforms in $\R^3$, $R_1\mu(x,y,z)$, defined as
$$R_1\mu(x,y,z)=\int \frac z{|(u,v,0)-(x,y,z)|^3}d\mu(u+iv).$$
This transform is the planar analogue of the Poisson transform. In particular, 
$$\lim_{z\to 0+} R_1\mu (x,y,z)=\frac {d\mu}{dm_2}(x+iy)$$
for all points $w=x+iy\in \C$ where  the Radon derivative
$$\frac {d\mu}{dm_2}(w)=\lim_{r\to 0+}\frac{\mu(B(w,r))}{|B(w,r)|}$$
exists.

For measures on the line or on the circle their Poisson integrals and Radon derivatives
(with respect to the one-dimensional Lebesgue measure) are very much related but not 
always equivalent. When the asymptotics of the Poisson integral and the ratio from the definition of the
Radon derivative are different near
a certain point it usually means that the measure is "irregular" near that point. 
It is not difficult to show that if $\mu$ is absolutely continuous then at a Lebesgue point
of its density function the Radon derivative of $\mu$ and the Poisson integral of $|\mu|$ (or $R_1|\mu|$ if $n>1$) behave equivalently.  Even for singular measures on the circle, if a measure possesses a certain symmetry near a point, then
the same equivalent behavior takes place, as follows for instance from \cite{AAN}, Lemma 4.1. In fact, it is 
not easy to construct a measure so that its Poisson integral and Radon derivative behaved
differently near a large set of points. The same can be said about the Riesz transform and the Radon derivative.
Thus
one may interpret our next result as an evidence that, for a planar measure $\mu$, most points where $C^\mu=0$  are "irregular."

\begin{theorem}\label{t4}
Let $\mu\in M$ and let $w=x+iy$ be a point of density (with respect to $m_2$) of the set $E=\{C^\mu=0\}$.
Then 
$$ \frac{\mu(B(w,r))}{\pi r^2}=o\left(R_1|\mu| (x,y,r)\right)\ \ as\ \ r\to 0+.$$
\end{theorem}

In view of the above discussion this implies

\begin{corollary} If $w$ is a point of density of the set $E=\{C^\mu=0\}$, such that there exists the Radon derivative $ {d|\mu|}/{dm_2}(w)\not = 0$, then
\begin{equation}\label{e11}\mu(B(w,r))=o\left(|\mu|(B(w,r))\right)\ \ as\ \ r\to 0+\end{equation}
and $ {d\mu}/{dm_2}(w)=0$.

\end{corollary}

Since $m_2$-almost every point of a set is its density point, 
we also obtain the following version of the result from \cite{TV}:

\begin{corollary}\label{V}
The set $E=\{C^\mu=0\}$ has measure zero with respect to the absolutely continuous component of $\mu$.
\end{corollary}

\begin{proof}[Proof of Theorem \ref{t4}]
without loss of generality $w=0$. Choose a $C_0^{\infty}$ test-function $\phi$ supported in $B:=B(0,r)$, and such that 
$0\leq \phi\leq D/r^2, |\nabla\phi|\leq A/r^3$ and $\int_{\C}\phi\,d m_2=1$. Denote the complement of $E$ by $E^c$.
Then

\begin{equation}\label{e12}
\int \phi d\mu = \langle \phi,\bar{\partial}C^\mu\rangle =\langle \bar{\partial}\phi,C^\mu\rangle=\langle \chi_{E^c}\bar{\partial}\phi,C^\mu\rangle
=\int \left(\int\frac{\chi_{E^c}\bar{\partial}\phi\,dm_2(z)}{\zeta -z}\right)d\mu(\zeta)
\end{equation}

All we need is to show that the last integral is small. Then, since the first integral in \eqref{e12} is similar to the right-hand side of \eqref{e11} we will
complete the proof.
The main idea for the rest of the proof is to make the function $F(\zeta)=\int\frac{\chi_{E^c}\bar{\partial}\phi\,dm_2(z)}{\zeta -z}$ "small" by subtracting
a linear combination of Cauchy kernels corresponding to points from $E$, which will not change its integral with respect to $\mu$. 

Namely,
let $a,b\in B(0,r)\cap E$ be any two points such that $|a-b|>r$. By the previous lemma for any $z\in B(0,r)$ there exist constants $A=A(z),B=B(z)$, of modulus at most 2, such that
\eqref{est} holds with $c=z$. Integrating \eqref{est} with respect to $\chi_{E^c}\bar{\partial}\phi\,dm_2(z)$ we obtain that 
$$\left|\int\frac{\chi_{E^c}\bar{\partial}\phi\,dm_2(z)}{\zeta -z}-\frac {A^*}{\zeta-a}-\frac {B^*}{\zeta-b}\right|<C\frac{\e(r)r}{|\zeta|^3}$$
outside of $B(0,2r)$ for some constants $A^*, B^*$, where $\e(r)=|B(0,r)\cap E^c|/r^2=o(1)$ as $r\to 0$.
The constants satisfy $|A^*|, |B^*|< 2\frac {\e(r)}r$.

Notice that if $w\in E$ then
$\int\frac 1{\zeta-w}d\mu=0 $ by the definition of the set $E$. Hence, since $a,b\in E$,

$$\int \left(\int\frac{\chi_{E^c}\bar{\partial}\phi\,dm_2(z)}{\zeta -z}\right)d\mu(\zeta)=
\int \left(\int\frac{\chi_{E^c}\bar{\partial}\phi\,dm_2(z)}{\zeta -z}-\frac {A^*}{\zeta-a}-\frac {B^*}{\zeta-b}\right)d\mu(\zeta)
$$$$
=\int_{B(0,2r)}+
\int_{\C\setminus B(0,2r)}=I_1+I_2.
$$

For $I_2$ we now have
$$\left|\int_{\C\setminus B(0,2r)} \left(\int\frac{\chi_{E^c}\bar{\partial}\phi\,dm_2(z)}{\zeta -z}-\frac {A^*}{\zeta-a}-\frac {B^*}{\zeta-b}\right)d\mu(\zeta)\right|
$$$$\leq C\int_{\C\setminus B(0,2r)}\frac{\e(r)r}{|\zeta|^3}d|\mu|(\zeta)\leq C\e(r)R_1|\mu|(0,0,r).
$$

In $I_1$ we estimate each summand separately. First,
$$\left|\int_{B(0,2r)}\left(\int\frac{\chi_{E^c}\bar{\partial}\phi\,dm_2(z)}{\zeta -z}\right)d\mu(\zeta)\right|\leq \int_{B(0,2r)}\frac D{r^3} \int\frac 1{|\zeta-z|}\chi_{E^c}dm_2(z)d|\mu|(\zeta)$$$$\leq C\frac{\sqrt{\e(r)}}{r^2}|\mu|(B(0,2r))\leq C\sqrt{\e(r)}R_1|\mu|(0,0,r).
$$
To estimate the second and third summands of $I_1$, recall that the only restriction on the choice of $a,b\in B(0,r)\cap E$ was that $|a-b|>r$. This condition
will be satisfied, for instance, if $a\in B_1=B(-\frac56 r,\frac 16r)$ and $b\in B_2=B(\frac56 r,\frac 16r)$. If we average the modulus of the second summand
over all choices of $a\in B_1\cap E$, recalling that $A^*=A^*(a)$ always satisfies $|A^*|\leq 2\frac{\e(r)}r$, we get
$$ \frac1{|B_1\cap E|}\int_{B_1\cap E}\left|\int_{B(0,2r)}\frac {A^*(a)}{\zeta-a}d\mu(\zeta)\right|dm_2(a)\leq
 \frac1{|B_1\cap E|}\int_{B(0,2r)}\int_{B_1\cap E}\frac {|A^*(a)|}{|\zeta-a|}dm_2(a)d|\mu|(\zeta)$$$$
\leq C \frac 1{r^2}\frac{\e(r)}{r}r|\mu|(B(0,2r))\leq
C\e(r)R_1|\mu|(0,0,r).
$$
It is left to choose $a\in B_1\cap E$ for which the modulus is no greater than its average. The same can be done for $b$.
The proof is finished.
\end{proof}

\section{Reflectionless measures and Combs}
\label{refl}

As was mentioned in the introduction, following \cite{BBEIM}, we will call a non-trivial continuous finite
measure $\mu\in M(\C)$ {\it reflectionless} if $C^{\mu}(z) = 0 $ at $\mu$-a.e. point $z$. 

Perhaps the simplest example of a reflectionless measure is the measure
$\mu=\frac 1\pi(1-x^2)^{-1/2}dx$ on $[-1,1]$, the harmonic measure of $\C\setminus [-1,1]$ corresponding
to infinity. The fact that $\mu$ is reflectionless can be verified through routine
calculations or via the conformal map interpretation of the harmonic measure.
It will also follow from a more general Theorem \ref{harmrefl1} below.

At the same time, since $C^\mu_*\asymp(1-x^2)^{-1/2}$ on $[-1,1]$, this simple example 
complements the statement of Theorem \ref{zero}.  Since the function $(1-x^2)^{-1/2}$
belongs to the "weak" $L^1(|\mu|)$, the summability condition for the Cauchy maximal function 
proves to be exact in its scale.

\vspace{.1in}

In the rest of this section we discuss further examples and properties of positive reflectionless measures
on the line.

\vspace{.1in}

Let us recall that functions holomorphic in the upper half plane $\C_+$  and mapping it to itself (having non-negative imaginary part) are called Nevanlinna functions. Let $M_+(\R)$ denote the class of finite positive measures compactly supported on $\R$. The function $f$  is a Nevanlinna function if and only if it has a form
$$
f(z) =az +b + \int_{\R}[\frac1{t-z} - \frac{t}{t^2+1}]d\rho(t)\,,
$$
where $\rho$ is a positive measure on $\R$ such that $\int \frac{d\rho(t)}{t^2+1}<\infty$, $a>0, b\in \R$ are constants. If the representing measure is from $M_+(\R)$ and  $f(\infty)=0$, the formula becomes simpler: $f(z)=\int\frac{d\mu(x)}{x-z}$.

\vspace{.2in}

\noindent {\bf Definition.} A simply connected  domain $\OO$ is {\it comb-like} if it is a subset of a half-strip \newline $ \{w:\Im w \in (0,\pi), \Re w>q\}$, 
for some $q\in\R$,
contains another half-strip
$\{w:\Im w \in (0,\pi), \Re w >r\}$ for some $r\in \R$
 and has the property that
\begin{equation}
\label{comblike}
\text{for any}\,\, w_0=u_0+iv_0\in \OO\,\,\text{the whole ray} \,\,\{w=u+iv_0, u\geq u_0\}\,\,\text{lies in}\,\, \OO\,.
\end{equation}
 If in addition $\HH^1(\partial\OO\cap B(0,R))<\infty$ for all finite $R$, we say that
$\OO$ is a {\it rectifiable comb-like} domain.

\vspace{.1in}

Let $\OO$ be a {\it rectifiable comb-like} domain, $\Gamma=\partial\OO$. Then by the Besicovitch theory  we know  that for $\HH^1$-a.e. pont $w\in \Gamma$ there exists an approximate tangent line to $\Gamma$, see \cite{B} for details. We wish to consider {\it rectifiable comb-like} domains satisfying the following geometric property:
\begin{equation}
\label{VH}
\text{for a.e.}\,\, w\in \Gamma\,\,\text{ approximate tangent line is either vertical or horizontal}.
\end{equation}

\vspace{.1in}

It is not difficult to verify that for any conformal map $F:\C_+\rightarrow \OO$, $\OO$ is comblike if and only if $F'$ is a Cauchy potential of $\mu\in M_+(\R)$: $F'(z)=\int\frac{d\mu(x)}{x-z}$. It is, therefore, natural to ask the following

\vspace{.1in}

\noindent{\bf Question}. Which {\it comb-like} domains correspond to reflectionless measures $\mu\in M_+(\R)$?

\vspace{.1in}

An answer would give a  geometric description of reflectionless measures from $M_+(\R)$. If, in addition, a comb-like domain is rectifiable, then the answer is given by

\begin{theorem}
\label{rectif}
1) Rectifiable comb-like domains
correspond exactly to those measures \newline $\mu\in M_+(\R)$ that  are absolutely continuous with respect to $dx$ and satisfy
\begin{equation}
\label{H1}
\int\frac{d\mu(x)}{x-z}\in H^1_{loc}(\C_+).
\end{equation}
2) An absolutely continuous measure satisfying \eqref{H1} is reflectionless if and only if
the corresponding comb-like domain has the property \eqref{VH}.
\end{theorem}

\vspace{.1in}

\noindent{\bf Remarks.} 

1) Of course not every comb-like domain gives rise to a reflectionless measure from $M_+(\R)$. Just take any comb-like domain which appears as $F(\C_+)$, where $F=\int^z\int\frac{d\mu(x)}{x-z}$ for a singular $\mu\in M_+(\R)$. By a result from \cite{JP} singular measures cannot be reflectionless.

2) On the other hand, even if $\mu=g(x)dx$ is a reflectionless absolutely continuous measure,
the corresponding conformal map $F=\int^z\int\frac{d\mu(x)}{x-z}:\C_+\rightarrow \OO$ can be
onto a non-rectifiable domain. 

3) For non-rectifiable domains we have no criteria to recognize which ones correspond to  reflectionless measures.

4) It is well known, and not difficult to prove, that the antiderivative of a  Nevanlinna function is a conformal map, see for instance \cite{Dur}. If $F=\int^z\int\frac{d\mu(x)}{x-z}, \mu\in M_+(\R)$ then $\Im F(x)$ is an increasing function on $\R$ whose derivative in the sense of distributions is $\mu$.
The image $F(\C_+)$ lies in the strip $\{\Im w \in (0,\pi\|\mu\|)\}$.

\vspace{.1in}

Theorem \ref{rectif} will follow from Theorems \ref{VHto} and \ref{toVH} below.

\begin{theorem}
\label{VHto}
Let $F$ be a conformal map of $\C_+$ on a rectifiable comb-like domain
 $\OO$. Then $F(z) =\int^z \int\frac{d\mu(x)}{x-z}$, $\mu\in M_+(\R), \mu<<dx$.
 Also $\int\frac{d\mu(x)}{x-z}\in H^1_{loc}(\C_+)$.
 If in addition $\OO$ satisfies \eqref{VH} then $\mu$ is reflectionless.
\end{theorem}

\vspace{.1in}

\begin{proof} without loss of generality $\OO\subset \{\Re z>0\}$. Put $\Phi = e^{F}$. Then the  image $\Phi(\OO)$ is  the subdomain of the  complement of the unit half-disk in $\C_+$ which is the union of rays $(R(\theta)e^{i\theta}, \infty)$. Consider the subdomain of the upper half-disk $D:=\{z: 1/z\in \Phi(\OO)\}$. Define $G$
as the smallest open domain containing $D$ and its reflection $\overline{D}:=\{\bar{z}: z \in D\}$. Then $G$ is a star-like domain inside the unit disk.
The preimage of $G\cap\R$ under $\Phi$ is the union of two  Infinite rays $R_1=[-\infty,a), R_2=(b,\infty], a<b$. Therefore, by reflection principle $\C\setminus [a,b]$ is mapped conformally (by the extension of $\Phi$ which we will also denote by $\Phi$) onto star-like $G$. 


\vspace{.1in}

 Since $\Phi: \C_+\rightarrow G$, where $G$ is star-like, it is well-known that $\arg \Phi(x+i\delta)$ is an increasing function of $x$, see \cite{Go}. 

\vspace{.1in}

We conclude that the argument of $\Phi$ is monotone. Therefore, $\Im F(x+i\delta)$ is monotone,
and so $\Im f(x+i\delta)$ is positive, where $f=F'$. We see that $f=F'$ is a Nevanlinna function. From the structure of our comb-like domain, we conclude immediately that its representing measure $\mu$ has compact support, so we are in $M_+(\R)$. Also, let us prove that $\mu<< dx$. The boundary of our comb is locally rectifiable. So $f=F'$ belongs locally to the Hardy class $H^1(\C_+)$, \cite{Privalov}. Since $\Im f$ is the Poisson integral of $\mu$,
$$\Im f =  P_{\mu}=\frac1{\pi}\int\frac{y}{(x-t)^2 +y^2} \,d\mu(t),$$
 and $f$ is in  $H^1(\C_+)$ locally, we conclude that $\mu =\Im f dx, \Im f \geq 0$ a.e., \cite{Privalov}. 

\vspace{.1in}

Now suppose that, in addition, $\OO=F(\C_+)$ has the property \eqref{VH}. Let us recall that for a simply connected domain with rectifiable boundary $\Gamma$ the restriction of the Hausdorff measure $\HH^1|_\Gamma$ is equivalent to the harmonic measure $\nu$ on $\OO$. Therefore the tangent lines to $\Gamma$ are either vertical or horizontal a.e. with respect to $\nu$.
The measure $\nu$ is the image of the harmonic measure $\lambda$ of $\C_+$ which is equivalent to the Lebesgue measure on the line. We have a conformal map $F$ (a continuous function up to the boundary of $\C_+$ because it is an anti-derivative of an $H^1_{loc}$-function) which pushes forward $\lambda$ to $\nu$. Call a point $w_0\in \Gamma$ accessible from $\OO$ if there exists a ray $x_0+iy, 0<y<1,$ such that
$w_0=\lim_{y\rightarrow 0} F(x_0+iy)$. Almost every point of $\Gamma$ (w.r. to $\nu$) is accessible from $\OO$. For $\nu$-a.e. accessible $w_0\in \Gamma$ where the tangent line is vertical (horizontal) we can say that $\Re F'(x_0) =0$ ($\Im F'(x) =0$). So $\R=E_1\cup E_2\cup E_3$, where $|E_3| =0, |E_1\cap E_2|=0$, and $E_1=\{x\in \R: \Re F'(x) =0\}$,
$E_2=\{x\in \R: \Im F'(x) =0\}$.  We already know that the measure $\mu=\Im F'(x) dx$ represents $f(z)=F'(z) = \int_{\R\setminus E_2}\frac{d\mu(t)}{t-z}$. Notice that $\int_{\R\setminus E_2}\cdot = \int_{E_1}\cdot$. But we also know that boundary values exist $dx$-almost everywhere, i.e. 
$$
\lim_{y\rightarrow 0}\Re \int_{E_1}\frac{d\mu(t)}{t-x-iy} =\Re F'(x) =0
$$
 for a.e. $x\in E_1$ and therefore for  $\mu$-a.e. $x\in E_1$. This means (see \cite{Privalov}) that
$$
p.v.\int_{\R}\frac{d\mu(x)}{x-z} =0\,\,\mu\text{-a.e.}
$$



\end{proof}

\noindent{\bf Definition.} A simply connected  {\it rectifiable comb-like} 
domain $\OO$ 
is called a comb if  its ``left" boundary consists of countably many horizontal and vertical segments.

A comb
is called a straight comb if $\OO = \{w:\Im w \in (0,\pi), \Re w >0\}\setminus S$, where the  set $S$ is relatively closed with respect to the strip $\{w:\Im w \in (0,\pi), \Re w >0\}$ and is the union of countably many horizontal intervals $R_n =(iy_n, l_n + i y_n]$.
We require also that
$$
\sum_n l_n <\infty\,.
$$

\vspace{.1in}

\noindent{\bf Example.}
Let $F$ be a conformal map of $\C_+$ on a comb $\OO$. By our last theorem $F'(z) = \int\frac{d\mu(x)}{x-z}$, where $\mu\in M_+(\R)$ is reflectionless: $C^{\mu}(x) = 0$ for $\mu$-a.e. $x$. 

\vspace{.1in}

\noindent{\bf Definition.} Let $E$ be  a compact subset of the real line. Let $E$ have positive logarithmic capacity, so Green's function $G$ of $\C\setminus E$ exists. The domain $\C\setminus E$ is called Widom domain if 
$$
\sum G(c) <\infty\,,
$$
where the summation goes over all critical points of $G$ (we assume that $G$ is a Green's function with pole at infinity.

\vspace{.1in}

\noindent{\bf Example.}
Let $E$ be  a compact subset of the real line of the positive length. We assume that
every point of $E$ is regular in the sense of Dirichlet for the domain $\C\setminus E$, and we also assume that $\C\setminus E$ is {\it not} a Widom domain. Such $E$ exist in abundance.
We will see below, that the harmonic measure $\omega$ of $\C\setminus E$ (with pole at infinity) is reflectionless. Consider $F(z) = \int^z\int\frac{d\omega(x)}{z-x}$ for $z\in \C_+$.
It is easy to see that $F(z) = G(z) + i \widetilde{G}(z) + const$, where $\widetilde{G}$ is the harmonic conjugate of $G$. This $F$ is a conformal map  (see \cite{Dur}) of $\C_=$ onto a domain $D$ lying in the strip $\{w: \Im w \in (0,\pi)\}$. It is easy to see that complementary intervals of $E$ will be mapped by $F$ onto straight horizontal segments on the boundary of $D$. Each finite complementary interval contains exactly one critical point of $G$, and clearly
the length of the corresponding straight horizontal segment is $G(c)$ (this follows from the formula $F(z) = G(z) + i \widetilde{G}(z) + const$). 

As the domain $\C\setminus E$ was {\it not} a Widom domain, we have that the sum of lengths
of abovementioned straight horizontal segment is infinite. So domain $D$ is not rectifiable.
Therefore the reflectionless property of $\mu$ alone does not say anything about the rectifiability of the domain, which is the target domain of the conformal map  $F(z) = \int^z\int\frac{d\mu(x)}{z-x}$.

\vspace{.3in}

\begin{theorem}
\label{toVH}
Let $\mu$ be absolutely continuous positive measure on $\R$ and let $C^{\mu}\in H^1_{loc}(\C_+)$. 
Then $F(z)=\int^z \int\frac{d\mu(x)}{x-z}$ is a conformal map of $\C_+$ onto a rectifiable comb-like domain $\OO$. If $\mu$ is reflectionless then $\OO$ has the property \eqref{VH}. 
\end{theorem}

\vspace{.2in}

\begin{proof}
Consider $F(z)=\int^z \int\frac{d\mu(x)}{x-z}$. Since $\mu$ is positive, it is a conformal map. If $\mu$ is such that $f(z)=C^{\mu}\in H^1_{loc}(\C_+)$ then $F(z)=\int^z f$ maps $\C_+$ onto a domain with locally rectifiable boundary (see \cite{Privalov}).

If, in addition, $\mu =\Im f dx$ is reflectionless, then for a.e. point of
$P:=\{x\in\R: \Im f(x)>0\}$ we have $\Re f(x) =0$. Conformal map $F(z)$ is continuous up to the boundary of $\C_+$ and its boundary values $F(x)$ form a (locally) absolutely continuous function, $F'(x) = f(x)$ a.e. As at almost every point we have either $\Im F'(x)=0$  or 
$\Re F'(x)=0$ we conclude that $\OO= F(C_+)$ has the property \eqref{VH}.

\end{proof}

\vspace{.1in}

\noindent We also need  the following definition.

\vspace{.1in}

\noindent{\bf Definition.} A compact subset $E$ in $\R$ is called homogeneous if there exist
$r,\delta>0$ such that for all $x\in E$, $|E\cap(x-h,x+h)|\geq \delta h$ for all $h\in (0,r)$.

\vspace{.1in}

\noindent{\bf Example.}
Let $E\subset\R$ be  a compact set of positive length. Let $\mu$ be a reflectionless measure supported on $E$, $\mu=g(x) dx$. Let in addition $E$ be a homogeneous set. Then $F(z) =\int^z \int\frac{d\mu(x)}{x-z}$ is a conformal map from $\C_+$ on a  rectifiable comb-like domain satisfying \eqref{VH}.

\vspace{.1in}

\begin{proof} The Cauchy integral $C^{gdx}$ considered in $\C\setminus E$ will be in the Hardy class $H^1(\C\setminus E)$. In fact the reflectionless property of $gdx$ implies that its limits from $C_\pm$ will be both integrable with respect to $dx|_E$. 

Now we use homogenuity of $E$ and Zinsmeister's theorem \cite{Z} to conclude that
$f(z)=C^{gdx}(z)$ is in the usual $H^1_{loc}(\C)$. Then the conformal map $F(z)=\int^z f$ maps $C_+$ onto a rectifiable subdomain of a strip. We use Theorem \ref{toVH} to get the rest of our example's claims. 
\end{proof}

\ms

The simple example of a reflectionless measure mentioned at the beginning of this section, as well as many other explicit examples,
are given by our next statement.

\begin{theorem}
\label{harmrefl1}
Let $E$ be  a compact set of positive lenght, $E\subset\R$.   Let $\omega$ be a harmonic measure of $\C\setminus E$ with pole at infinity. 
Then $\omega$ is reflectionless.
\end{theorem}

\bs

\noindent{\bf Example.} The simplest comb is a strip $\{w:\Im w \in (0,\pi), \Re w >0\}$.
Consider $F(z) =\log (z+\sqrt{z^2-1})$. It maps conformally $\C_+$ onto the strip. Its derivative $f(z) = \frac1{\sqrt{z^2-1}}$ is $\frac1{\pi}\int\frac{dx}{\sqrt{1-x^2}}\frac1{x-z}$ and $d\mu=\frac1{\pi}\frac{dx}{\sqrt{1-x^2}}$ is the harmonic measure of $\C\setminus [-1,1]$.

\vspace{.2in}

\begin{proof}[Proof of Theorem \ref{harmrefl1}] 
We need to show that $C^\omega=0$ at $\omega$-a.e point. From our definitions it can be seen, that $C^\omega$ on the line coincides with the Hilbert transform
of $\omega$, which in its turn is asymptotically equivalent to the conjugate Poisson transform $Q^\omega$. Thus all we need to establish is that

\begin{equation}
Q^\omega(x+ih)=\int_{\R}\frac{x-y}{(x-y)^2 +h^2}d\omega(y)=\Re \int\frac{d\omega(y)}{x-ih-y}\to 0\ \ {\text as }\ \ h\to 0+
\end{equation}
for almost every $x$.
Instead, we have that the Green's function $F(x)$ defined as
$$F(x)=\int\log |x-y|d\omega(y) + C_\infty,$$
where $C_\infty$ is a real constant (Robin's constant), is equal to 0 at every density point of $E$, see for example \cite{Hay}.
The idea of the proof is to show that $Q^\omega(x+i\e)$ behaves like $(F(x+\e)+F(x-\e))/\e$ near almost every $x$. The technical details
are as follows.

Introduce
\begin{equation}
\label{phi}
\phi (y) := \frac12 \log\frac{|1-y|}{|1+y|} + \frac{y}{y^2+1}\,,
\end{equation}
$$
\phi_{x, h}(y) := \frac1h\phi(\frac{y-x}h)\,.
$$
The function $\phi(y)$ decreases as $1/y^2$ at infinity, hence it is in $L^1(\R,dx)$ and so are $\phi_{x, h}(y)$ with a uniform bound on the norm. However, these functions 
are not bounded, which makes it difficult to use them in our estimates. To finish the proof we will first obtain a bounded version of $\phi_{x, h}(y)$ through the following
averaging procedure.

Let $\omega=g(x)dx$. Choose  $x$ to be a Lebesgue point of $g$ and a density point of $E$. Fixing sufficiently small $h>0$ we can find the set $A(x,h)\subset (x-h,x-h/2)\cup (x+h/2, x+h)$ such that 
\begin{itemize}
\item $ A(x,h)$ consists of density points of $E$,
\item $|A(x,h)| \geq h/2$,
\item $A(x,h)$ is symmetric with respect to $x$.
\end{itemize}
Let $T_{x,h}:=T:=\{t\in (0,h): x+t \in A(x,h)\}$. Then $|T| \geq h/4$.
Now put
$$
\psi_{x,h}(y) := \frac1{|T|}\int_T \phi_{x,t} (y)\,dt\,.
$$
By  \eqref{phi} one can see immediately that

\begin{equation}
\label{psi}
|\psi_{x,h}|\leq \frac{M}{h}\ \ \text{ for some } \ M>0\ \text{ and }\ \ |\psi_{x,h}(y)|\leq C\frac{h}{y^2},\ \  \text{ for }\ |y|>h\ .
\end{equation}

 Also, since
$$
\int \phi\,dy =0\,.
$$
we have that 
$$
\int \psi_{x,h}\,dy =0\,.
$$
Therefore,
$$
|\int g(y) \psi_{x,h}(y)\,dy|=|\int (g(y)-g(x)) \psi_{x,h}(y)\,dy|\leq \int |g(y)-g(x)|| \psi_{x,h}|(y)\,dy. $$

Now notice that \eqref{psi} implies that $| \psi_{x,h}|$ is majorated by an approximate unity (for instance, by a constant multiple of the Poisson kernel
corresponding to $z=x+ih$). Since $x$ is a Lebesgue point for $g(x)$, this means that the last integral tends to 0 as $h\to0$.

 Looking  at the definitions of $T_{x,h}$ and
$\psi_{x,h}(y)$ we can see that
$$
\int_{\R} g(y)\psi_{x,h}(y)\,dy= \frac1{|T_{x,h}|}\int_{T_{x,h}}\bigg[\frac1{2t}(F(x+t) -F(x-t)) - \Re \int\frac{g(y)dy}{x-it-y}\bigg]\,,
$$
where $F(x)$ is the Green's function. As we mentioned before, $F$ is zero at the density points of  $E$. We conclude that 
$$
\Re \frac1{|T_{x,h}|}\int_{T_{x,h}}dt\int\frac{g(y)dy}{x-it-y}\rightarrow 0,\,\,\, h\rightarrow 0+\,.
$$
for a.e. $x$ on the Borel support of $g$. Since the Cauchy integral of $g$ has a limit a.e. we obtain that

$$
\Re\int\frac{g(y)dy}{x-ih-y}\rightarrow 0,\,\,\, h\rightarrow 0+.
$$
\end{proof}

%
%
%

\noindent{\bf Remark.}  All reflectionless measures on $\R$ discussed in this section, including
those provided by Theorem \ref{harmrefl1} are absolutely continuous
with respect to Lebesgue measure. One may wonder if there exist singular reflectionless measures. The answer is negative. More generally, as follows from a theorem from \cite{JP}, if principal values of the Hilbert transform exist $\mu$-a.e. for 
a continuous $\mu\in M(\R)$ then $\mu<<dx$ .

\vspace{.1in}

\indent


\begin{thebibliography}{99}

\bibitem{AAN} A. B. Alexandrov, J. M. Anderson, A. Nicolau. Inner functions, Bloch spaces and symmetric measures,
Proc. London Math. Soc. (3) 79 (1999), no. 2, 318--352.


\bibitem{BBEIM} 
E.D. Belokolos, A.I. Bobenko, V.Z. Enol'skii, A.R. Its and V.B. Matveev,
Algebro-Geometric Approach to Nonlinear Integrable Equations,  
Springer, Berlin (1994).

\bibitem{B} P. Mattila, Geometry of sets and measures in Euclidean spaces. Fractals and rectifiability. Cambridge Studies in Advanced Mathematics, 44. Cambridge University Press, Cambridge, 1995. xii+343 pp.

\bibitem{Dur} P. L. Duren, Univalent Functions, Grundlehren der Mathematischen Wissenschaften [Fundamental Principles of Mathematical Sciences], 259. Springer-Verlag, New York, 1983.

\bibitem{EG} L.C. Evans, R. Gariepy, Measure Theory and Fine Properties of Functions,
Studies in Advanced Mathematics. CRC Press, Boca Raton, FL, 1992.

\bibitem{Gam}
T. W. Gamelin, Uniform Algebras, Prentice-Hall, Inc., Englewood Cliffs, N. J., 1969.

 
\bibitem{Go} G. M. Golusin, Geometric theory of functions. Hochschulbcher fr Mathematik, Bd. 31. VEB Deutscher Verlag der Wissenschaften, Berlin, 1957. xii+438 pp. 30.0X 


\bibitem{Hay} W. K. Hayman, P. B. Kennedy, Subharmonic Functions, vol. 1, Academic Press, 
London-New York, 1976.

\bibitem{JP} P. Jones, A. Poltoratski, Asymptotic growth of Cauchy transforms, Ann. Acad. Sci. Fenn. Math, 2004 


 \bibitem{KN} M. Krein, A. Nudelman
{\em The Markov moment problem and extremal problems. Ideas and
problems of P. L. \v Ceby\v sev and A. A. Markov and their further
development}. Translated from the Russian by D. Louvish.
Translations of Mathematical Monographs, Vol. 50. American
Mathematical Society, Providence, R.I., 1977. v+417 pp.

\bibitem{MV} M. Melnikov, J. Verdera, A geometric proof of the $L\sp 2$ boundedness of the Cauchy integral on Lipschitz graphs,  Internat. Math. Res. Notices  1995,  no. 7, 325--331.

\bibitem{NTV}F. Nazarov, S. Treil, A. Volberg , Cauchy integral and Calder�-Zygmund operators on nonhomogeneous spaces,
Int. Math. Res. Not. 15 (1997) 703726.

\bibitem{T} X. Tolsa, $L\sp 2$ -boundedness of the Cauchy integral operator for continuous measures, Duke Math. J. 98 (1999), 269-304.

\bibitem{TV} X. Tolsa, J. Verdera, May the Cauchy transform of a non-trivial finite measure vanish on the support of the measure?  Ann. Acad. Sci. Fenn. Math.  31  (2006),  no. 2, 479--494.

\bibitem{Z} M. Zinsmeister, Espaces de Hardy et domaines de Denjoy. (French) [Hardy spaces and Denjoy domains] Ark. Mat. 27 (1989), no. 2, 363--378.

\bibitem{Privalov} I. Privalov, Boundary properties of analytic functions,
Gosudarstv. Izdat. Tehn.-Teor. Lit., Moscow-Leningrad, 1950. 336 pp.

\end{thebibliography}
\end{document}